\documentclass[a4paper,11pt]{article}
\usepackage{latexsym}
\usepackage{amsmath}
\usepackage{amssymb}
\usepackage{graphicx}
\usepackage[round]{natbib}
\usepackage{delarray}
\usepackage[usenames]{color}
\usepackage{colortbl}
\newtheorem{proposition}{Proposition}[section]

\newtheorem{thm}{Theorem}[section]
\newcommand\sign{\mathop{\mathrm{sign}}}
\title{Integral representations of one dimensional projections for 
multivariate stable densities}
\author{Muneya MATSUI${}^{\dag}$ and Akimichi TAKEMURA${}^{\ddag 
}$ \\ \\
$\dag$ 
\small Department of Mathematics, Keio University \\
$\ddag$ \small Graduate School of Information Science and Technology,
University of Tokyo
}
\date{August, 2006}
\begin{document}
\maketitle
\begin{abstract}
  We consider the numerical evaluation of one dimensional projections of
  general multivariate stable densities introduced by
  \cite{Abdul-Hamid:Nolan:1998}.  In their approach higher order
  derivatives of one dimensional densities are used, which seem
  to be cumbersome in practice. Furthermore there are some
  difficulties for even dimensions.  In order to overcome these
  difficulties we obtain the explicit finite-interval integral
  representation of one dimensional projections for all dimensions.
  For this purpose we utilize the imaginary part of complex
  integration, 
  whose real part corresponds to the derivative of the one dimensional
  inversion formula.
  We also give summaries on relations between various
  parameterizations of stable multivariate density and its one
  dimensional projection.
\end{abstract}

\section{Introduction}
Stable distributions are known as the limiting distributions of the
general central limit theorem and the time invariant distributions for L\'evy
processes. Because of these probabilistically very important properties, many
studies have been carried out on their theoretical
aspects. Comprehensive books have been published 
(\cite{Araujo:Gine:1980} or \cite{Christoph:Wolf:1993}).
Furthermore 
applications of stable distributions 
to heavy-tailed data have been growing over past few decades.
These applications have appeared in many fields like finance, Internet
traffic or physics. Especially in time series, the observations from the
causal stable ARMA or fractional
stable ARIMA model can be considered as a multivariate stable random vector.
We refer to \cite{Uchaikin:Zolotarev:1999},
\cite{Rachev:2003} or \cite{Adler:Feldman:Taqqu:1998} for these applications.
As various multivariate heavy tailed data have become available in these fields,
the treatments of multivariate
stable distributions are needed.
However, the multivariate stable 
distributions require handling of their spectral measures which is computationally difficult.

We briefly review some important recent progresses concerning the  
multivariate stable distributions. 
Concerning the calculation of the integral with respect to the
spectral measure, \cite{Byczkowski:Nolan:Rajput:1993} have
approximated the spectral measure by a discrete spectral
measure which is uniformly close to the original spectral measure. 
They also have given several instructive graphs
of the stable densities.
Concerning the approximations, 
\cite{Davydov:Nagaev:2002} have outlined possible directions although
their results are still theoretical. For the estimation of the parameters, 
see \cite{Marcus:Luis:2003} and \cite{Davydov:Paulauskas:1999} and 
for the hypothesis testing see \cite{Mittnik:Rachev:Ruschendorf:1999}.
Although there are many other researches, 
accurate calculations of the multivariate stable densities are clearly of
basic importance.

For the calculations of the multivariate stable densities, 
the method of \cite{Abdul-Hamid:Nolan:1998}
based on  one dimensional projection is very promising.
In the approach of \cite{Abdul-Hamid:Nolan:1998}, the basic ingredient
is a function $g_{\alpha,d}(v,\beta)$ defined by 
\begin{equation}
\label{eq:g-alpha-d}
g_{\alpha,d}(v,\beta)
= 
\begin{cases}
\frac{1}{(2\pi)^d}\int_0^\infty \cos
\left(vu-\left(\beta\tan\frac{\pi\alpha}{2}\right)
u^\alpha\right)u^{d-1}e^{-u^\alpha}du, & \alpha \neq 1, \\
\frac{1}{(2\pi)^d}\int_0^\infty \cos
 \left(vu+\frac{2}{\pi}\beta u\log u\right)u^{d-1}e^{-u}du, & \alpha=1,
\end{cases}
\end{equation}
where $\alpha$ is the characteristic exponent of the stable
distribution, $d$ is the dimension and $\beta$ corresponds to the
skewness parameter.

Evaluating $g_{\alpha,d}(v,\beta)$ using this very definition 
is not very practical because the integrand is
infinitely oscillating and changing its sign.  
Therefore an alternative evaluation of 
$g_{\alpha,d}(v,\beta)$ is desirable.

In this paper, we obtain explicit finite-interval integral
representations of $g_{\alpha,d}(v,\beta)$, which do not require the 
derivatives of one dimensional densities.  
The integrand in our representation is well-behaved without infinite
oscillations.  Furthermore, our integral representation covers 
all dimensions and all parameter values, except for an exceptional
case of $\alpha \le 1, \beta =1$ and even dimension, where we need an
additional term which is easily tractable.

The calculation of $g_{\alpha,d}(v,\beta)$ corresponds to
the evaluation of the equation (2.2.18)
of \cite{Zolotarev:1986}, which is stated without an explicit proof.
For odd dimension, the real part of the equation is equivalent to $g_{\alpha,d}(v,\beta)$
and for even dimension, the imaginary part is equivalent to $g_{\alpha,d}(v,\beta)$.
Complete derivations of the equation (2.2.18) of \cite{Zolotarev:1986}
is also useful for the numerical evaluations of the derivatives of one dimensional
symmetric densities investigated in \cite{Matsui:Takemura:2006}.

This paper is organized as follows. For the rest of this section we
prepare preliminary results based on \cite{Abdul-Hamid:Nolan:1998}
and \cite{Nolan:1998}.
In Section \ref{sec:finite-integral} we derive 
finite-interval integral representations of one dimensional
projections.
In Section \ref{sec:various-projections} we summarize 
relations between
various parameterizations of the multivariate stable density and its 
one dimensional projection as propositions.
Some discussions and directions for further studies are given 
in Section \ref{sec:procedures-future-work}.
Some proofs are postponed to Section \ref{sec:proofs}. 

\subsection{Definitions and preliminary results}
\label{sec:definitions}
Let $\mathbf{X} =(X_1,X_2,\ldots,X_d)$ denote a $d$-dimensional  $\alpha$-stable random
vector with the characteristic function 
$$
\Phi_{\alpha,d}(\mathbf{t}) = E \exp (i\langle \mathbf{t}, \mathbf{x} \rangle ),
$$ 
where $\mathbf{t}=(t_1,t_2,\ldots,t_d)$ and $\langle\ ,\ \rangle$ denotes
the usual inner product.
Among several definitions of stable distributions, we give the spectral representation of
$\Phi_{\alpha,d}(\mathbf{t})$.
This  requires integration over the unit sphere
$\mathbb{S}^d=\{\mathbf{s} : \Vert\mathbf{s}\Vert=1\}$ 
in $\mathbb{R}^d$ with respect to a spectral measure
$\Gamma$ 
(Theorem 2.3.1 of \cite{Samorodnitsky:Taqqu:1994}, Theorem 14.10 of \cite{Sato:1999}).
The spectral definition is given as
\begin{equation}
\label{eq:definition-characteristic-a-representation}
\Phi_{\alpha,d}(\mathbf{t})= 
\begin{cases}
\exp\left(-\int_{\mathbb{S}^d} | \langle
				  \mathbf{t},\mathbf{s} \rangle |^\alpha
				  \left(1-i\sign(\langle
				   \mathbf{t}, \mathbf{s} \rangle)
				  \tan\frac{\pi\alpha}{2}\right)\Gamma(d\mathbf{s})+i\langle
				 \mathbf{t}, \mathbf{\mu} \rangle
                               \right), & \alpha \neq 1, \\
\exp\left(-\int_{\mathbb{S}^d} | \langle
				  \mathbf{t},\mathbf{s} \rangle |
				  \left(1+i\frac{2}{\pi}\sign(\langle
				   \mathbf{t}, \mathbf{s} \rangle)\log|\langle
				   \mathbf{t},\mathbf{s}\rangle| \right)\Gamma(d\mathbf{s})+i\langle
				 \mathbf{t}, \mathbf{\mu} \rangle
                               \right), & \alpha=1.
\end{cases}
\end{equation}
The pair $(\Gamma,\mu)$ is unique. 
The definition (\ref{eq:definition-characteristic-a-representation})
corresponds to Zolotarev's (A) representation and 
we use the notation
$S_{\alpha,d}(\Gamma,\mu)$ for
this definition following \cite{Nolan:1998}.

For another definition, we give a multivariate version of Zolotarev's
(M) parameterization, which is also defined in \cite{Nolan:1998}. We use
notation $S^M_{\alpha,d}(\Gamma,\mu_0)$. 
\begin{align}
\label{eq:definition-characteristic-m-representation}
& \Phi_{\alpha,d}(\mathbf{t}) \nonumber \\
& \quad = \begin{cases}
\exp\left(-\int_{\mathbb{S}^d} | \langle
				  \mathbf{t},\mathbf{s} \rangle |^\alpha
				  \left(1+i\sign(\langle
				   \mathbf{t}, \mathbf{s}\rangle) \tan\frac{\pi\alpha}{2}( |\langle
				   \mathbf{t},\mathbf{s} \rangle|^{1-\alpha}-1) 
				  \right)\Gamma(d\mathbf{s})+i\langle
				 \mathbf{t}, \mathbf{\mu_0} \rangle
                               \right), 
& \alpha\neq 1, \\
\exp\left(-\int_{\mathbb{S}^d} | \langle
				  \mathbf{t},\mathbf{s} \rangle |
				  \left(1+i\frac{2}{\pi}\sign(\langle
				   \mathbf{t}, \mathbf{s} \rangle)\log|\langle
				   \mathbf{t},\mathbf{s}\rangle| \right)\Gamma(d\mathbf{s})+i\langle
				 \mathbf{t}, \mathbf{\mu_0} \rangle
                               \right), 
& \alpha=1.
\end{cases}
\end{align}
In the equation (\ref{eq:definition-characteristic-m-representation}),
we avoid the discontinuity at $\alpha=1$ and the
divergence of the 
mode as $\alpha \to 1$ which are observed in the definition 
(\ref{eq:definition-characteristic-a-representation})
in the non-symmetric case. The difference of (M)
parameterization from 
(A) representation is only 
$$
\mu_0=\mu-\tan\frac{\pi\alpha}{2}\int_{\mathbb{S}^d}\mathbf{s}\
\Gamma(d\mathbf{s}), \qquad \alpha\neq 1.
$$

Now consider inverting the characteristic function
$\Phi_{\alpha,d}(\mathbf{t})$ to evaluate the multivariate stable density.
The useful idea of \cite{Nolan:1998} and \cite{Abdul-Hamid:Nolan:1998} is to perform this
inversion with respect to the polar coordinates.
The projection of an $\alpha$-stable random vector onto one dimensional
stable random variable is explained as follows. 
Here we follow the arguments of Ex.2.3.4 of \cite{Samorodnitsky:Taqqu:1994}.
For any $\mathbf{t}\in \mathbb{R}^d$ consider the characteristic
function $E\exp(iu\langle \mathbf{t},\mathbf{x}\rangle)$ of 
a linear combination $Y=\sum_i t_i X_i$,
i.e.,
\begin{eqnarray*}
&&E\exp(iu\langle \mathbf{t},\mathbf{x}-\mu \rangle) \\  
&& \quad 
=
\begin{cases}
\exp\left(-|u|^\alpha\left[ \int_{\mathbb{S}^d} | \langle
				  \mathbf{t},\mathbf{s} \rangle |^\alpha\Gamma(d\mathbf{s})
			  -i\sign(u)\tan\frac{\pi\alpha}{2}
\int_{\mathbb{S}^d}\sign(\langle
			   \mathbf{t}, \mathbf{s} \rangle)| \langle
				  \mathbf{t},\mathbf{s} \rangle |^\alpha
			  \Gamma(d\mathbf{s})\right]  \right),
&
\alpha \neq 1, \\
\exp\left(-|u| \left[\int_{\mathbb{S}^d} | \langle
			  \mathbf{t},\mathbf{s} \rangle | \Gamma(d\mathbf{s})
                          +i\frac{2}{\pi}\sign(u)\int_{\mathbb{S}^d} \langle
			  \mathbf{t},\mathbf{s} \rangle \left(\log|u|+\log|\langle
			   \mathbf{t},\mathbf{s}\rangle|\right)\Gamma(d\mathbf{s})\right]
    \right), 
&\alpha = 1.
\end{cases}
\end{eqnarray*}
Then, from the definition below, $E\exp(iu\langle
\mathbf{t},\mathbf{x}\rangle)$ is considered as a characteristic function
of one dimensional stable distribution with the following parameters.
\begin{align}
\label{def-A-sigma}
\sigma(\mathbf{t})&=\left(\int_{\mathbb{S}^d}|\langle \mathbf{t}, \mathbf{s}
 \rangle|^\alpha \Gamma (d \mathbf{s}) \right)^{1/\alpha},\\
\label{def-A-beta}
\beta(\mathbf{t})&=(\sigma(\mathbf{t}))^{-\alpha}\int_{\mathbb{S}^d} \sign
\langle \mathbf{t}, \mathbf{s} \rangle\ |\langle \mathbf{t}, \mathbf{s}
 \rangle|^\alpha\ \Gamma (d \mathbf{s}), \\
\label{def-A-mu}
\mu(\mathbf{t}) &=
\left\{
\begin{array}{ll} 
0,& \alpha \neq 1, \\
-\frac{2}{\pi}\int_{\mathbb{S}^d}\langle \mathbf{t}, \mathbf{s}
 \rangle \log| \langle \mathbf{t}, \mathbf{s}
 \rangle|\ \Gamma (d \mathbf{s}), & \alpha=1.
\end{array}
\right.
\end{align}
Regarding $u$ as the length and $\mathbf{t}$ as the angle, we 
transform the $d$-dimensional rectangular integral into
the polar coordinate integral. In the following theorem of
\cite{Nolan:1998}, multivariate stable densities of (A) representation
are projected onto  one dimensional densities in (A)
representation. Since in the original theorem (\cite{Nolan:1998})
there are some trivial typographical errors,
we correct them in the formula (\ref{eq:g-alpha-d}) and the following theorem.
\begin{thm} {\rm (\cite{Nolan:1998})} \quad
\label{thm:A-A-representation}
Let $0<\alpha<2$ and $\mathbf{X} =(X_1,X_2,\ldots,X_d)$ be
an $\alpha$-stable random vector $S_{\alpha,d}(\Gamma,\nu)$ with $d\ge
 2$. Then the density $f_{\alpha,d}(\mathbf{x})$ of $\mathbf{X}$ is given as 
\begin{equation*}
f_{\alpha,d}(\mathbf{x})= 
\left\{
\begin{array}{ll}
\int_{\mathbb{S}^d} g_{\alpha,d}\left(
\displaystyle
\frac{ \langle \mathbf{x}-\nu,\mathbf{s} \rangle
}{\sigma(\mathbf{s})},\beta(\mathbf{s})\right)\left(\sigma(\mathbf{s})\right)^{-d} \ d \mathbf{s},
 & \alpha\neq 1,\\
\int_{\mathbb{S}^d}g_{1,\alpha}\left(
\displaystyle
\frac{\langle \mathbf{x}-\nu,\mathbf{s} \rangle -\mu(\mathbf{s})-(2/\pi)\beta(\mathbf{s})\sigma(\mathbf{s})\log\sigma(\mathbf{s})}{\sigma(\mathbf{s})},\beta(\mathbf{s})
\right)\left(\sigma(\mathbf{s})\right)^{-d} \ d \mathbf{s},
 & \alpha=1,
\end{array}\right.
\end{equation*}
where $g_{\alpha,d}(v,\beta)$ is given in (\ref{eq:g-alpha-d}).
\end{thm}

\cite{Abdul-Hamid:Nolan:1998} expressed
$g_{\alpha,d}(v,\beta)$ as an integral over a finite interval, utilizing
the integral expressions of one dimensional densities. 
They explained the method when $d$ is odd as follows. 
{}From Theorem 2.2.3 of \cite{Zolotarev:1986}, the one dimensional density $g_{\alpha,1}(v,\beta)$
can be written as
\begin{equation*}
g_{\alpha,1}(v,\beta)=\int_a^b h_g(\theta;\alpha,\beta,v) d\theta,
\end{equation*}
where $-\pi/2 \le a<b \le \pi/2$ and $h_g$ is a somewhat complicated
but explicit function. Then, for odd $d$, differentiating
$g_{\alpha,1}$  $d-1$ times  with respect to $v$ gives the representation
$$
g_{\alpha,d}(v,\beta)= c(d)\int_a^b
\frac{\partial^{d-1}h_g(\theta;\alpha,\beta,v)}{\partial v^{d-1}}d \theta,
$$
where $c(d)$ is some constant.  
However, the computation of
$\partial^{d-1}h_g(\theta;\alpha,\beta,v)/\partial v^{d-1}$ seems to
be cumbersome. Furthermore, they did not give the corresponding
expression for even $d$.
In the following, we obtain explicit finite-interval integral
representations of $g_{\alpha,d}(v,\beta)$, which do not require the
derivatives of one dimensional densities.  


\section{Finite-interval integral representations of projected 
one dimensional functions}
\label{sec:finite-integral}
In this section we obtain the finite-interval integral representations of 
$g_{\alpha,d}(v,\beta)$ in (\ref{eq:g-alpha-d}).
It is based on the function
$h^n(x;\alpha,\beta)$ defined below in 
(\ref{eq:n-derivative-invesion}), which 
corresponds to the $n$-th order derivative of one
dimensional inversion formula including the imaginary part.

We first prepare some notations which are the same as in
\cite{Zolotarev:1986}.
\allowdisplaybreaks{
\begin{align*}
K(\alpha)&=\alpha-1+\sign (1-\alpha), \qquad \alpha\neq 1,\\
\theta&=\beta K(\alpha)/\alpha, \qquad \alpha\neq 1,\\
U_\alpha(\varphi;\theta)&=\left(\frac{\sin\frac{\pi}{2}(\varphi+\theta)}{\cos\frac{\pi}{2}\varphi}\right)^{\alpha/(1-\alpha)}
\frac{\cos\frac{\pi}{2}\{(\alpha-1)\varphi+\alpha\theta\}}{\cos\frac{\pi}{2}\varphi},
\qquad \alpha \neq 1,\\
U_1(\varphi;\beta)&= \frac{\pi}{2} \frac{1+\beta\varphi}{\cos\frac{\pi}{2}}
\exp\left( \frac{\pi}{2} \left(\varphi+\frac{1}{\beta} \right)
  \tan\frac{\pi}{2}\varphi\right), \\
r(\varphi) &= 
\left\{
\begin{array}{ll}
\left(
\frac{\sin\alpha(\varphi+\pi\theta/2)}{x\cos\varphi}
\right)^{1/(1-\alpha)},
&\alpha \neq 1, \\ 
\exp(-x/\beta+(\varphi+\pi/(2\beta))\tan \varphi),
&\alpha =1,
\end{array}
\right. \\
\tau &= 
\left\{
\begin{array}{ll}
(\alpha/x)^{1/(1-\alpha)},
&\alpha \neq 1, \\ 
\exp(-x-1),
&\alpha =1,
\end{array}
\right. \\
V_n(\varphi) &= r^n\left(\frac{\pi}{2}\varphi\right)
\left\{ r'\left(\frac{\pi}{2}\varphi
	  \right)\sin\frac{\pi}{2}(n+1)(\varphi+1)+r\left(\frac{\pi}{2}\varphi\right)\cos\frac{\pi}{2}(n+1)(\varphi+1) \right\},\\
W_n(\varphi) &= r^n\left(\frac{\pi}{2}\varphi\right)
\left\{ r\left(\frac{\pi}{2}\varphi \right)\sin\frac{\pi}{2}(n+1)(\varphi+1)-r'\left(\frac{\pi}{2}\varphi\right)\cos\frac{\pi}{2}(n+1)(\varphi+1) \right\}.
\end{align*}

We now define a function $h^n(x;\alpha,\beta)$  as  
\begin{eqnarray}
\label{eq:n-derivative-invesion}
h^n(x;\alpha,\beta) &=& \frac{1}{\pi}\int_0^\infty (iz)^n
 \exp\left(izx+\psi(z;\alpha,-\beta)\right)dz, 
\end{eqnarray} 
where 
\begin{eqnarray*}
\psi(z;\alpha,\beta)  &=&
\left\{
\begin{array}{ll}
-z^{\alpha}\exp(-i\frac{\pi}{2}\theta\alpha),
&\alpha \neq 1, \\
-\frac{\pi}{2}-i\beta z \log z,
&\alpha =1,
\end{array}
\right.
\end{eqnarray*}
is the logarithm of the characteristic function of
stable distributions which are analytically extended from the semi-axis
$\mathrm{Re}\;z>0$ into the complex plane. Then,
$g_{\alpha,d}(v,\beta)$  in (\ref{eq:g-alpha-d}) is expressed in terms
of $h^n(x;\alpha,\beta)$ as follows. Note that $\beta$ in
$g_{\alpha,d}(v,\beta)$ is not equal to $\beta$ in
$h^n(x;\alpha,\beta)$, because 
$\beta$ in $h^n(x;\alpha,\beta)$ is based on one dimensional (B)
parameterization. Therefore, we need the following definitions,
$$
\beta_B = \frac{2}{\pi K(\alpha)} \arctan \left(\beta \tan \frac{\pi \alpha}{2} \right),\quad
\theta_B = \beta_B K(\alpha)/\alpha, \qquad \alpha\neq 1.\\
$$
\begin{proposition}
\label{lem:relation-g-h}
Define
 $x=\left(\cos(\pi/2\alpha\theta_B)\right)^{1/\alpha}v$ and
 $y=\pi/2v+\beta\log (\pi/2)$. Then a
  representation of $g_{\alpha,d}(v,\beta)$  in (\ref{eq:g-alpha-d})
in terms of $h(\cdot\ ;\alpha,\beta)$ is given as follows:
\begin{eqnarray*}
g_{\alpha,d}(v,\beta) &=&
\left\{
\begin{array}{ll} 
\mathrm{Re}\ 
 \frac{\left(\cos\left(\pi/2\alpha\theta_B
		 \right)\right)^{d/\alpha}}{2^d(\pi i)^{d-1}} h^{d-1}(|x|;\alpha,\beta_B^\ast),
 & \alpha \neq 1,\\ &
\\
\mathrm{Re}\ 
 \frac{\pi}{4^d i^{d-1}} h^{d-1}(y^\ast;\alpha,|\beta|),
 &  \alpha=1,
\end{array}
\right.
\end{eqnarray*}
where $\beta_B^\ast=\beta_B\ \sign x$ and $y^\ast=y\ \sign \beta$.
\end{proposition}

Proof of this proposition is given in Section \ref{sec:proofs}.
Although the imaginary number $i$ is present in the above formulas, it
is for notational convenience and we do not have to write
$\mathrm{Re}h(\cdot)$ or $\mathrm{Im}h(\cdot)$ depending on the
dimension $d$.

By Proposition \ref{lem:relation-g-h}, it suffices to express 
$h^n(x;\alpha,\beta)$ as a finite-interval integral.
We now give our main theorem.
We mention that the following results are
partially stated in \cite{Zolotarev:1986} for the real case without
proofs and the equation (2.2.34) of \cite{Zolotarev:1986} is not complete.
Therefore we give the complete results including the calculation of
imaginary parts. The full proof is given in Section \ref{sec:proofs}.

\begin{thm}
\label{thm:differential-density}
Let $x>0$ in the case $\alpha \neq 1$ and $\beta>0$ in the case
 $\alpha=1$. 
Then the finite-interval integral representations of 
 $h^n(x;\alpha,\beta)$ are as follows. \\
{\rm (a)}\ $(\alpha \neq 1,\beta \neq 1)$ or $(\alpha>1, \beta=1)$
\begin{eqnarray*}
\mathrm{Re}\ h^n(x;\alpha,\beta) &=& \frac{1}{2}\int_{-\theta}^1 \exp
 \left(-x^{\alpha/(\alpha-1)}U_\alpha(\varphi;\theta)\right)
 V_n(\varphi)d\varphi.  \\
\mathrm{Im}\ h^n(x;\alpha,\beta) &=& \frac{1}{2}\int_{-\theta}^1 \exp
 \left(-x^{\alpha/(\alpha-1)}U_\alpha(\varphi;\theta)\right)
 W_n(\varphi) d\varphi.
\end{eqnarray*}
{\rm (b)}\ $(\alpha=1,\beta\neq 1)$
\begin{eqnarray*}
\mathrm{Re}\ h^n(x;1,\beta) &=& \frac{1}{2}\int_{-1}^1 \exp
 \left(-e^{-x/\beta}U_1(\varphi;\beta)\right)
V_n(\varphi)d\varphi.  \\
\mathrm{Im}\ h^n(x;1,\beta) &=& \frac{1}{2}\int_{-1}^1 \exp
 \left(-e^{-x/\beta}U_1(\varphi;\beta)\right)
W_n(\varphi)d\varphi.  
\end{eqnarray*}
{\rm (c)}\ $(\alpha<1,\beta=1)$ \\
$\mathrm{Re}\ h^n(x;\alpha,1)$ is the same as in {\rm (a)}.
\begin{eqnarray*}
\mathrm{Im}\ h^n(x;\alpha,1) &=& \frac{1}{2}\int_{-\theta}^1 \exp
 \left(-x^{\alpha/(\alpha-1)}U_\alpha(\varphi;\theta)\right)
 W_n(\varphi) d\varphi -\frac{1}{\pi}\int_0^\tau \exp(xr-r^\alpha) r^n dr.
\end{eqnarray*}
{\rm (d)}\ $(\alpha=1,\beta=1)$ \\
$\mathrm{Re}\ h^n(x;1,1)$ is the same as in {\rm (b)}.
\begin{eqnarray*}
\mathrm{Im}\ h^n(x;1,1) &=& \frac{1}{2}\int_{-1}^1 \exp
 \left(-e^{-x/\beta}U_1(\varphi;\beta)\right) W_n(\varphi) d\varphi-\frac{1}{\pi}\int_0^\tau \exp\left(xr+r\log r\right)r^n dr.
\end{eqnarray*}
\end{thm}
Interestingly, the case $\beta=1$ is different from other values of
the parameters. As suggested by \cite{Zolotarev:1986}, the representations of
Theorem \ref{thm:differential-density} are essentially different from
the $n$-th order derivative  of
the integral representations of densities given in the equation (2.2.18)
of \cite{Zolotarev:1986}.  

\section{Relations between representations of
multivariate stable densities and one dimensional projections}
\label{sec:various-projections}
In this section, we present propositions concerning the
projections of a stable random vector onto one dimensional stable random
variable. 
The idea of projection is useful for the purpose of the computation of
the multivariate densities.
There exist two representations (A) and (M) for the
multivariate stable distributions,  whereas three
representations (A), (B) and (M) are known for general one dimensional
stable distributions. The projection of the multivariate representation
(A) onto one dimensional representation (A) is given in Theorem
\ref{thm:A-A-representation}, which is obtained by \cite{Nolan:1998}.  
Although \cite{Nolan:1998} mentioned the projection of (M) onto (M)
as a theorem without proof, his is different from ours. Therefore we describe the projection of
the multivariate representation (M) onto the projected one dimensional representation (M)
with a precise proof. Note that the idea of projection entirely belongs
to \cite{Nolan:1998}. Further we present results for all projections without
proofs by adding the other 4 projections.
Note that for 
$\alpha=1$  multivariate representations (A) and (M) and one dimensional
representations (A) and (M) are the same.  Therefore we omit
statements
on the projections (A $\to$ M), (M $\to$ A) and (M $\to$ B) for $\alpha=1$.

\subsection{One-dimensional projections in various representations}
Here we present propositions on various projections. Note that the
parameters $\sigma(\mathbf{t})$, $\beta(\mathbf{t})$ and
$\mu(\mathbf{t})$ used in this section have already been defined by (\ref{def-A-sigma}),
(\ref{def-A-beta}) and (\ref{def-A-mu}) in Section \ref{sec:definitions}, respectively.
\begin{proposition}[The projection (M $\to$ M)]
\label{thm:M-M-representation}
Let $0<\alpha<2$ and $\mathbf{X} =(X_1,X_2,\ldots,X_d)$ be an $\alpha$-stable random vector $S_{\alpha,d}^M(\Gamma,\nu)$ with $d\ge 2$. Then,
 the density $f_{\alpha,d}^M(\mathbf{x})$ is given as follows. 
\begin{eqnarray*}
f_{\alpha,d}^M( \mathbf{x} )
=\int_{\mathbb{S}^d}g_{\alpha,d}^M\left(\frac{ \langle \mathbf{x}-\nu,\mathbf{s} \rangle
-\mu_M^M(\mathbf{s})}{\sigma(\mathbf{s})},\beta(\mathbf{s})\right)
\left(\sigma(\mathbf{s})\right)^{-d} \ d \mathbf{s},
\end{eqnarray*}
where \\
\begin{eqnarray*}
g_{\alpha,d}^M(v,\beta)=\frac{1}{(2\pi)^d}\int_0^\infty \cos \left(vu+
\left(\beta\tan\frac{\pi\alpha}{2}(u-u^\alpha)\right)\right)u^{d-1}e^{-u^\alpha}du
\end{eqnarray*}
and
$$
\mu_M^M(\mathbf{t})=\tan\frac{\pi\alpha}{2}\left(
\beta(\mathbf{t}) \sigma(\mathbf{t})-
\int_{\mathbb{S}^d} \langle \mathbf{t}, \mathbf{s}
 \rangle \ \Gamma (d \mathbf{s})\right).
$$ \end{proposition}

When $\alpha=1$, (A) parametrization in Theorem
\ref{thm:A-A-representation} and (M)
parametrization in Proposition \ref{thm:M-M-representation} are  the same.
It is easy to confirm $f_{1,d}^M(
\mathbf{x})=f_{1,d}( \mathbf{x})$ and
$g_{1,d}^M(v,\beta)=g_{1,d}(v,\beta)$ only by notational change.
Note that at $\alpha=1$, $\mu_M^M(\mathbf{t})$ becomes
\begin{eqnarray*}
\label{eq:mu-at-alpha=1}
\mu_M^M(\mathbf{t})&=&\frac{2}{\pi}\left(
\beta(\mathbf{t}) \sigma(\mathbf{t})\log \sigma(\mathbf{t})-
\int_{\mathbb{S}^d} \langle \mathbf{t}, \mathbf{s}
 \rangle \ \log | \langle \mathbf{t}, \mathbf{s}
 \rangle |\Gamma (d \mathbf{s})  \right) \\
&=&
\mu(\mathbf{s})+(2/\pi)\beta(\mathbf{s})\sigma(\mathbf{s})\log\sigma(\mathbf{s}). \nonumber
\end{eqnarray*}
\begin{proposition}[The projection (A $\to$ B)]
\label{thm:A-B-representation}
Under the same conditions and notations as in Theorem 
\ref{thm:A-A-representation}, the density $f_{\alpha,d}( \mathbf{x} )$
 is given as follows.
\\
\begin{equation*}
f_{\alpha,d}(\mathbf{x})= 
\left\{
\begin{array}{ll}
\int_{\mathbb{S}^d}g_{\alpha,d}^B\left(\frac{ \langle \mathbf{x}-\nu,\mathbf{s} \rangle
}{\sigma_B(\mathbf{s})},\beta_B(\mathbf{s})\right)
\left(\sigma_B(\mathbf{s})\right)^{-d} \ d \mathbf{s},
& \alpha\neq 1,\\
\int_{\mathbb{S}^d}g_{1,d}^B\left(\frac{ \langle \mathbf{x}-\nu,\mathbf{s} \rangle
- \mu_B(\mathbf{s})}{\sigma_B(\mathbf{s})},\beta(\mathbf{s})\right)
\left(\sigma_B(\mathbf{s})\right)^{-d}\ d \mathbf{s},
 & \alpha=1,
\end{array}\right.
\end{equation*}
where
\begin{equation*}
g_{\alpha,d}^B(v,\beta)
= 
\begin{cases}
\frac{1}{(2\pi)^d}\int_0^\infty \cos \left(vu-
u^\alpha \sin \left(\frac{\pi}{2}K(\alpha)\beta \right)
\right)u^{d-1}e^{-u^\alpha\cos\left(\pi/2K(\alpha)\beta_B\right)}du,
& \alpha \neq 1, \\
\frac{1}{(2\pi)^d}\int_0^\infty \cos \left(vu+ \beta u \log u \right)u^{d-1}e^{-\pi/2u}du,
 & \alpha=1,
\end{cases}
\end{equation*}
and
\begin{equation*}
\sigma_B(\mathbf{t})
= 
\begin{cases}
\frac{ \sigma(\mathbf{t}) }{ \left( \cos\left(
 \frac{\pi}{2} K(\alpha) \beta_B(\mathbf{t}) \right)
 \right)^{1/\alpha}},
& \alpha \neq 1, \\
\frac{2}{\pi}\sigma(\mathbf{t}),
 & \alpha=1,
\end{cases}
\end{equation*}
$$
\beta_B(\mathbf{t}) = \frac{2}{\pi K(\alpha)} \arctan \left(\beta (
 \mathbf{t}) \tan
 \frac{\pi \alpha}{2} \right),
$$
$$\mu_B(\mathbf{t})=\sigma_B(\mathbf{t})\beta(\mathbf{t})\log \sigma_B(\mathbf{t})-\frac{2}{\pi}\int_{\mathbb{S}^d} \langle \mathbf{t}, \mathbf{s}
 \rangle \ \log | \langle \mathbf{t}, \mathbf{s}
 \rangle |\Gamma (d \mathbf{s}).$$ 
\end{proposition}
\begin{proposition}[The projection (A $\to$ M)]
Under the same conditions and relations as in Theorem
 \ref{thm:A-A-representation}, the density $f_{\alpha,d}(\mathbf{x})$ is
 given as follows. 
\begin{eqnarray*}
f_{\alpha,d}( \mathbf{x} )
=\int_{\mathbb{S}^d}g_{\alpha,d}^M\left(\frac{ \langle \mathbf{x}-\nu,\mathbf{s} \rangle
-\mu_M(\mathbf{s})}{\sigma(\mathbf{s})},\beta(\mathbf{s})\right)
\left(\sigma(\mathbf{s})\right)^{-d} \ d \mathbf{s}, \quad 
\alpha \neq 1,
\end{eqnarray*}
where
$$\mu_M(\mathbf{t})=\sigma(\mathbf{t})\beta(\mathbf{t})
\tan \frac{\pi\alpha}{2}. $$
\end{proposition}
\begin{proposition}[The projection (M $\to$ A)]
Under the same conditions and relations as in Proposition
 \ref{thm:M-M-representation}, the density $f^M_{\alpha,d}(\mathbf{x})$ is
 given as follows. 
\begin{eqnarray*}
f_{\alpha,d}^M( \mathbf{x} )
=\int_{\mathbb{S}^d}g_{\alpha,d}\left(\frac{ \langle \mathbf{x}-\nu,\mathbf{s} \rangle
-\mu^M_A(\mathbf{s})}{\sigma(\mathbf{s})},\beta(\mathbf{s})\right)
\left(\sigma(\mathbf{s})\right)^{-d} \ d \mathbf{s}, \quad 
\alpha \neq 1,
\end{eqnarray*}
where
$$\mu_A^M(\mathbf{t})=-\tan\frac{\pi\alpha}{2} \int_{\mathbb{S}^d} \langle \mathbf{t}, \mathbf{s}
 \rangle \ \Gamma (d \mathbf{s}). $$
\end{proposition}
\begin{proposition}[The projection (M $\to$ B)]
Under the same conditions and relations as in Proposition \ref{thm:M-M-representation}, the density $f^M_{\alpha,d}(\mathbf{x})$ is
 given as follows. 
\begin{eqnarray*}
f_{\alpha,d}^M( \mathbf{x} )
=\int_{\mathbb{S}^d}g^B_{\alpha,d}\left(\frac{ \langle \mathbf{x}-\nu,\mathbf{s} \rangle
-\mu^M_A(\mathbf{s})}{\sigma_B(\mathbf{s})},\beta_B(\mathbf{s})\right)
\left(\sigma_B(\mathbf{s})\right)^{-d} \ d \mathbf{s},\quad 
\alpha \neq 1.
\end{eqnarray*}
\end{proposition}

\subsection{Finite-integral representations for other representations}
In \cite{Zolotarev:1986} only the analytic extension of the
characteristic function defined by (B) representation is considered
and Theorem \ref{thm:differential-density} corresponds to (B) representation of stable
distributions. Therefore, we have introduced the parameters $\beta_B$
and $\theta_B$ in the finite integral-representation for (A) of
Proportion \ref{lem:relation-g-h}.
Although it is possible to derive (A) or (M) representation versions of
Theorem \ref{thm:differential-density} utilizing the results in \cite{Nolan:1997},
we have to consider the analytic extension of the characteristic function
defined by (A) or (M) representation which needs very complicated arguments.
Therefore, in this paper, we confine our results to the case of (B)
representation and for (A) and (M) representation, we utilize
Proposition \ref{lem:relation-g-h} and \ref{lem:relation-gM-h}. Note that, accordingly,
for $g_{\alpha,d}^B(v,\beta)$ in Proposition \ref{lem:relation-gB-h}, we do not need
to use the extra parameters like $\beta_B$ or $\theta_B$ other than that
defined in Theorem \ref{thm:differential-density}.

\begin{proposition}
\label{lem:relation-gM-h}
Let $g_{\alpha,d}^M(v,\beta)$ be as in Proposition
 \ref{thm:M-M-representation}. Define 
$x= \left(\cos\left(\pi/2\alpha\theta_B\right)\right)^{1/\alpha}\left(v+\tan\left(\pi/2\alpha\theta_B\right) \right)$ and $y$ as in
 Proposition \ref{lem:relation-g-h}. Then the  
  representation of $g^M_{\alpha,d}(v,\beta)$ using function
 $h(\cdot\ ;\alpha,\beta)$ is as follows.
\begin{eqnarray*}
g^M_{\alpha,d}(v,\beta) &=&
\left\{
\begin{array}{ll} 
\mathrm{Re}\ 
 \frac{\left(\cos\left(\pi/2\alpha\theta_B
		 \right)\right)^{d/\alpha}}{2^d(\pi i)^{d-1}} h^{d-1}(|x|;\alpha,\beta_B^\ast),

 &  \alpha \neq 1,\\ &
\\
g_{1,d}(v,\beta)
, &  \alpha=1,
\end{array}
\right.
\end{eqnarray*}
where $\beta_B^\ast=\beta\ \sign x$.
\end{proposition}

\begin{proposition}
\label{lem:relation-gB-h}
Let $g_{\alpha,d}^B(v,\beta)$ be as in Proposition
 \ref{thm:A-B-representation}. Then the representation of $g^B_{\alpha,d}(v,\beta)$ using function
 $h(\cdot\ ;\alpha,\beta)$ is as follows.
\begin{eqnarray*}
g^B_{\alpha,d}(v,\beta) &=&
\left\{
\begin{array}{ll} 
\mathrm{Re}\ 
 \frac{1}{2^d(\pi i)^{d-1}} h^{d-1}(|v|;\alpha,\beta^\ast),
 &  \alpha \neq 1,\\ &
\\
\mathrm{Re}\ 
 \frac{1}{2^d (\pi i)^{d-1}} h^{d-1}(v^\ast;1,|\beta|)
, &  \alpha =1,
\end{array}
\right.
\end{eqnarray*}
where $\beta^\ast=\beta\ \sign v$ and $v^\ast = v\
 \sign \beta$.
\end{proposition}

\section{Some discussions and future works}
\label{sec:procedures-future-work}
In this paper we focused on one dimensional projections $g$ of the general
multivariate stable density.  This is only one step in calculating the
density itself.
We need to substitute parameters which are functions on the unit
sphere $\mathbf{s}\in \mathbb{S}^d$ into $g$ and then
we need to integrate over the unit sphere
to evaluate the density itself. 
We presented improvements in evaluation of $g$.  Improvements 
of other steps are also important for the evaluation of the
multivariate density.
By further careful examinations of the integrand of the
function $h^n(x;\alpha,\beta)$, we might find some useful
regularities like one dimensional finite integral representations
stated in Section 3 of \cite{Nolan:1997}.

For our future work we consider showing
various properties like tail dependencies or relations between the spectral
measure and tails. Theoretically the asymptotic estimates of
multivariate stable densities are obtained by \cite{Watanabe:2000} or
\cite{Hiraba:2003}. Relations between these
theoretical results and our representations or numerical results are of our next concern.

Furthermore, we can also consider the expansions or
the asymptotic expansions of functions $g$'s for the boundary values. 
Many expansions of densities are found in \cite{Zolotarev:1986}. The method of
the expansions may be directly applied to the projected functions $g$'s.

\section{Proofs}
\label{sec:proofs}

In this section, we give the proofs of Proposition \ref{lem:relation-g-h}, 
Theorem \ref{thm:differential-density} and Theorem \ref{thm:M-M-representation}.

\subsection{Proof of Proposition \ref{lem:relation-g-h}}
For $\alpha \neq 1$, define
\begin{equation*}
Ig_{\alpha,d}(v,\beta)=\frac{1}{(2\pi)^d} \int_0^\infty \sin 
\left(vu-\left(\beta \tan\frac{\pi\alpha}{2}\right) u^\alpha
\right) u^{d-1} e^{-u^\alpha}du.
\end{equation*}
Then by simple notational change we have
\begin{eqnarray*}
g_{\alpha,d}(v,\beta)+iIg_{\alpha,d}(v,\beta)  
& =&\frac{1}{(2\pi)^d}\int_0^\infty
 e^{iuv-i\left( \beta\tan\frac{\pi\alpha}{2}\right) u^\alpha}u^{d-1}e^{-u^\alpha}
 du\\
& =& \frac{1}{2^d (\pi i)^{d-1}}\frac{1}{\pi}\int_0^\infty
(iu)^{d-1}e^{iuv}e^{-u^\alpha
\left\{1+i\tan\left(\frac{\pi}{2}\alpha \theta_B \right) \right \} } du\\
& =& \frac{1}{2^d (\pi i)^{d-1}}\frac{1}{\pi}\int_0^\infty
(iu)^{d-1}e^{iuv}
e^{-u^\alpha/\cos\left(\frac{\pi}{2}\alpha\theta_B\right) e^{i(\pi/2\alpha\theta_B)}}du.
\end{eqnarray*}
Replacing $u$ by
$u=\left(\cos\left(\frac{\pi}{2}\alpha\theta_B\right)\right)^{1/\alpha}t$, we have
\begin{eqnarray*}
&& \frac{1}{2^d(\pi i)^{d-1}}\left(\cos\left(\frac{\pi}{2}\alpha
		     \theta_B\right)\right)^{\alpha/d}\frac{1}{\pi}\int_0^\infty
(it)^{d-1}
e^{it \left( \cos\left(\frac{\pi}{2}\theta_B\right) \right)^{d/\alpha}v}
e^{-t^\alpha e^{i(\pi/2\alpha\theta_B)}}dt \\
&&\hspace{2mm} = \frac{1}{2^d(\pi i)^{d-1}}\left(\cos\left(\frac{\pi}{2}\alpha
		     \theta_B\right)\right)^{\alpha/d}\frac{1}{\pi} \int_0^\infty
(it)^{d-1} \exp (itx+\psi(t;\alpha,-\beta_B) ) dt \\
&&\hspace{2mm} = \frac{1}{2^d(\pi i)^{d-1}}\left(\cos\left(\frac{\pi}{2}\alpha
		     \theta_B\right)\right)^{\alpha/d} h^{d-1}(x;\alpha,\beta_B).
\end{eqnarray*}
Taking the real part of this equation, we obtain the desired result when
$x>0$. When $x<0$, replacing $v$ by $-v$ and $\beta$ by $-\beta$, we can
obtain the result in a similar fashion. Note that $g_{\alpha,d}(-v,-\beta)=g_{\alpha,d}(v,\beta)$.

For $\alpha=1$, define
$$
Ig_{1,d}(v,\beta)=\frac{1}{(2\pi)^d} \int_0^\infty
\sin\left(vu+\frac{2}{\pi}\beta u\log u\right)u^{d-1}e^{-u}du.
$$
What we need is to calculate 
$
g_{1,d} (v,\beta)+iIg_{1,d}(v,\beta)
$
and to consider the real part of this equation. Since the proof is similar
to the proof of the case $\alpha\neq 1$, 
we omit the rest of the proof for $\alpha = 1$.
 \hfill $\Box$

\subsection{Proof of Theorem \ref{thm:differential-density}}
The function $h^0(x;\alpha,\beta)$, which can be regarded as the
inversion formula including the imaginary part, has the following
representation ((2.2.20) of \cite{Zolotarev:1986}):
\begin{equation*}
h^0(x;\alpha,\beta)=\frac{1}{\pi}\int_0^\infty \exp\left(izx+\psi(z;\alpha,-\beta)\right)dz. 
\end{equation*}
Differentiating $h^0(x;\alpha,\beta)$ $n$ times with respect to $x$, we have
\begin{equation*}
h^n(x;\alpha,\beta)=\frac{1}{\pi}\int_0^\infty (iz)^n
\exp\left(izx+\psi(z;\alpha,-\beta)\right)dz. 
\end{equation*}
For the calculation of $h^n(x;\alpha,\beta)$ we consider the same contour as in \cite{Zolotarev:1986}, which is 
\begin{equation*}
 \Gamma =\left\{z:\mathrm{Im} \left(izx+\psi(z;\alpha,-\beta)\right)=0,\quad
  \frac{\pi}{2}k \le \arg z \le \frac{\pi}{2} \right\},
\end{equation*}
where 
\begin{eqnarray*}
k &=&
\left\{
\begin{array}{ccc}
-\theta        & \mbox{if}\quad \alpha \neq 1, &  \\ 
-1             & \mbox{if}\quad  \alpha =1,     &\beta \neq -1,   \\
1              & \mbox{if}\quad  \alpha \le 1,  &\beta = -1.
\end{array}
\right.
\end{eqnarray*}
Since readers can refer to
\cite{Zolotarev:1986} if necessary, the details of the contour are omitted.
{}From Lemma 2.2.3 of \cite{Zolotarev:1986}, we only have to calculate the integration along the contour $\Gamma$. 

\medskip
\noindent
(a)\ $(\alpha \neq 1,\beta \neq 1)$ or  $(\alpha>1, \beta=1)$. Direct
calculation gives
\begin{eqnarray}
h^n(x;\alpha,\beta)&=&\frac{1}{\pi}\mathrm{Re}\ \int_\Gamma
 \exp\left(izx+\psi(z;\alpha,-\beta)\right)(iz)^n dz \nonumber \\
&&+\frac{1}{\pi}\mathrm{Im} \int_\Gamma 
\exp\left(izx+\psi(z;\alpha,-\beta)\right)(iz)^n dz \nonumber \\
&=&\frac{1}{\pi}
\int_\Gamma \exp \left\{\mathrm{Re}\ (izx+\psi(z;\alpha,-\beta))\right\}\mathrm{Re}\ 
\{(iz)^n dz\} \nonumber \\
&&+\frac{1}{\pi}
\int_\Gamma \exp \left\{\mathrm{Re}\  (izx+\psi(z;\alpha,-\beta))\right\}\mathrm{Im}\ 
\{(iz)^n dz\} \nonumber \\
&=&
\frac{1}{\pi}\int_\Gamma \exp (-W(\varphi))\left(\mathrm{Re}\ \{(iz)^n dz \}+\mathrm{Im}\ 
\{(iz)^n dz\} \right), \label{eq:h-proof}
\end{eqnarray}
where
\begin{eqnarray*}
W(\varphi) &=&
\left\{
\begin{array}{ll}
x^{\alpha/(\alpha-1)} U_\alpha(2\varphi/\pi;\theta),
&\alpha \neq 1, \\ 
\exp(-x/\beta)U_1(2\varphi/\pi;\beta),
&\alpha =1.
\end{array}
\right.
\end{eqnarray*}
For details for derivation of $W(\varphi)$, see p.76 of \cite{Zolotarev:1986}.
Replacing $z$ by $z=re^{i\varphi}$ gives 
\begin{eqnarray*}
(iz)^n&=&(re^{i(\varphi+\frac{\pi}{2})})^n \\
&=& r^n e^{in(\varphi+\frac{\pi}{2})} \\
&=& r^n (\cos n(\varphi+\frac{\pi}{2})+i\sin n(\varphi+\frac{\pi}{2}) ).
\end{eqnarray*}
Since 
\begin{eqnarray*}
d(r\cdot \sin \varphi) &=& r'\sin\varphi \;d\varphi + r\cos\varphi \;d\varphi, \\
d(r\cdot \cos \varphi) &=& r'\cos\varphi  \;d\varphi - r\sin\varphi \;d\varphi,
\end{eqnarray*}
$$
dz = r'\cos \varphi  \;d\varphi 
 - r\sin\varphi  \;d\varphi +i\{r'\sin \varphi+r\cos \varphi\}  \;d\varphi.
$$
Combining the equations above, we have
\begin{eqnarray*}
\mathrm{Re}\ (iz)^n dz &=& r^n \left( r'\sin(n+1)\left(\varphi+\frac{\pi}{2}\right)+
		    r\cos(n+1)\left(\varphi+\frac{\pi}{2}\right)\right)d\varphi,
\\
\mathrm{Im}\ (iz)^n dz &=& r^n \left( r\sin(n+1)\left(\varphi+\frac{\pi}{2}\right)-
		    r'\cos(n+1)\left(\varphi+\frac{\pi}{2}\right)\right)d\varphi.
\end{eqnarray*}
{}From these equations and the relations of angles and the contour $\Gamma$
on p.76 of \cite{Zolotarev:1986}, 
($\ref{eq:h-proof}$) becomes
\begin{eqnarray*}
\mathrm{Re}\ h(x;\alpha,\beta) &=& \frac{1}{\pi}
 \int_{-\frac{\pi}{2}\theta}^\frac{\pi}{2} 
\exp\left(-x^{\alpha/(\alpha-1)}U_\alpha(2\varphi/\pi;\theta)\right)\\
&&\times r^n \left( r'\sin(n+1)\left(\varphi+\frac{\pi}{2}\right)+
		    r\cos(n+1)\left(\varphi+\frac{\pi}{2}\right)\right)
d\varphi \\
&=& \frac{1}{2}
 \int_{-\theta}^1 
\exp\left(-x^{\alpha/(\alpha-1)}U_\alpha(\varphi;\theta)\right)V_n(\varphi)
d\varphi
\end{eqnarray*}
and
\begin{eqnarray*}
\mathrm{Im}\ h(x;\alpha,\beta) &=& \frac{1}{\pi}
 \int_{-\frac{\pi}{2}\theta}^\frac{\pi}{2} 
\exp\left(-x^{\alpha/(\alpha-1)}U_\alpha(2\varphi/\pi;\theta)\right)\\
&& \times r^n \left( r\sin(n+1)\left(\varphi+\frac{\pi}{2}\right)-
		    r'\cos(n+1)\left(\varphi+\frac{\pi}{2}\right)\right)d\varphi \\
&=& \frac{1}{2}
 \int_{-\theta}^1 
\exp\left(-x^{\alpha/(\alpha-1)}U_\alpha(\varphi,\theta)\right)W_n(\varphi)
d\varphi.
\end{eqnarray*}

\medskip\noindent
(b)\ $(\alpha=1,\beta\neq 1)$. \ 
Since similar results hold in this case, we omit the
poof. 

\medskip\noindent
(c)\ $(\alpha<1,\beta=1)$. \ Here we must consider additional contour
$\Gamma^\ast=\{z:\mathrm{Re}\ z=0,-\tau \le \mathrm{Im}\ z \le 0\}$ as
stated in p.76 of \cite{Zolotarev:1986}. The contour
$\Gamma^\ast$ satisfies
$$
\mathrm{Im}\ ((izx)+\psi(z;\alpha,-1))=0
$$
for $z \in \Gamma^\ast$. The integral along $\Gamma^\ast$ becomes
$$
\frac{1}{\pi} \int_{\Gamma^\ast} \exp(izx+\psi(z;\alpha,-1))(iz)^n dz = 
\frac{1}{\pi} \int_{\Gamma^\ast} \exp(\mathrm{Re}\ (izx+\psi(z;\alpha,-1)))(iz)^n dz. 
$$
Replacing $z$ by $z=r e^{-\frac{\pi}{2}i}$ gives 
$$
\frac{1}{\pi}\int_0^\tau \exp\left(\mathrm{Re} (izx+\psi(z;\alpha,-1)) \right)
 r^n (-i) dr 
= -\frac{i}{\pi}\int_0^\tau \exp\left(-W\left(-\frac{\pi}{2}\right) \right)
 r^n dr,
$$
where
\begin{eqnarray*}
W\left(-\frac{\pi}{2}\right) &=&xr \sin\left(-\frac{\pi}{2}\right)+r^\alpha
 \cos \alpha \left(-\frac{\pi}{2}+\frac{\pi}{2}\theta \right) \\
&=& -xr+r^\alpha
\end{eqnarray*}
for $\alpha \neq 1$ and 
\begin{eqnarray*}
W\left(-\frac{\pi}{2}\right) &=& xr \sin
 \left(-\frac{\pi}{2}\right)+r\log r\sin
 \left(-\frac{\pi}{2}\right)+\left(-\frac{\pi}{2}+\frac{\pi}{2}r\cos
  \left(-\frac{\pi}{2}\right) \right) \\
&=& -xr-r\log r
\end{eqnarray*}
for $\alpha =1$. \hfill $\Box$

\subsection{Proof of Proposition \ref{thm:M-M-representation} }
Without loss of generality, we assume $\nu=0$ throughout the proof.
The characteristic function $\phi_{\alpha,d}^M(t)$ can be written
as 
\begin{eqnarray*}
\Phi_{\alpha,d}^M(\mathbf{t})&=&\exp\Bigg( -\int_{\mathbb{S}^d}|\langle \mathbf{t}, \mathbf{s}
 \rangle|^\alpha \Gamma (d \mathbf{s}) +
i \tan\frac{\pi\alpha}{2} \int_{\mathbb{S}^d} \sign
\langle \mathbf{t}, \mathbf{s} \rangle\ |\langle \mathbf{t}, \mathbf{s}
 \rangle|^\alpha\ \Gamma (d \mathbf{s}) \\
&&\hspace{1cm} 
-i\tan\frac{\pi\alpha}{2} \int_{\mathbb{S}^d} \sign
\langle \mathbf{t}, \mathbf{s} \rangle\ | \langle \mathbf{t}, \mathbf{s}
 \rangle |\ \Gamma (d \mathbf{s}) 
\Bigg)\\
&=&
\exp\left(-(\sigma(\mathbf{t}))^\alpha+i(\sigma(\mathbf{t}))^\alpha\tan
     \frac{\pi\alpha}{2}\beta(\mathbf{t})
-i\tan\frac{\pi\alpha}{2} \int_{\mathbb{S}^d}\langle \mathbf{t}, \mathbf{s}
 \rangle \ \Gamma (d \mathbf{s}) \right)\\
&=&
\exp\left(-(\sigma(\mathbf{t}))^\alpha-i(\sigma(\mathbf{t}))^\alpha\tan
     \frac{\pi\alpha}{2}\beta(\mathbf{t})\left((\sigma(\mathbf{t}))^{1-\alpha}-1 \right)+i\mu_M^M(\mathbf{t}) \right) \\
&=&\exp\left( -\psi^M\left(\mathbf{t}
\right) \right).
\end{eqnarray*}
Then, the inversion formula gives
\begin{eqnarray*}
f_{\alpha,d}^M(\mathbf{x})&=&
(2\pi)^{-d} 
\int_{\mathbb{R}^d }e^{-i\langle \mathbf{x},\mathbf{t} \rangle
}\exp\left( -\psi^M\left(\mathbf{t}
\right)\right) d \mathbf{t}. 
\end{eqnarray*}
Note that for a positive real number $r>0$ and a vector $\mathbf{s}\in
\mathbb{R}^d$, $\mu_M^M(r\mathbf{s})=r\mu_M^M(\mathbf{s})$,
$\sigma(r\mathbf{s})=r\sigma(\mathbf{s})$ and $\beta(r\mathbf{s})=\beta(\mathbf{s})$.
Putting $\mathbf{t}=r\mathbf{s}$ where $r>0$ and $\mathbf{s} \in
\mathbb{S}^d$, we obtain 
\begin{eqnarray*}
f_{\alpha,d}^M(\mathbf{x})
&=&(2\pi)^{-d}\int_{\mathbb{S}^d}\int_0^\infty
\exp\Big[-i\langle \mathbf{x},\mathbf{s}\rangle r-(r\sigma(\mathbf{s}))^\alpha-i(r\sigma(\mathbf{s}))^\alpha
\tan\frac{\pi\alpha}{2}\beta(\mathbf{s})\left( (r\sigma(\mathbf{s}))^{1-\alpha}-1 \right)\\
&& \hspace{3.8cm}
+ir\mu_M^M(\mathbf{s})\Big]
r^{d-1} dr d \mathbf{s}.
\end{eqnarray*}
Furthermore, replacing $r$ by $\mu=r \sigma(\mathbf{s})$, we get
\begin{eqnarray*}
f_{\alpha,d}^M(\mathbf{x})&=&(2\pi)^{-d}\int_{\mathbb{S}^d}\int_0^\infty
\exp\bigg[ -i\frac{\langle \mathbf{x},\mathbf{s}\rangle
-\mu_M^M(\mathbf{s})}{\sigma(\mathbf{s})} u
-\left\{u^\alpha+iu^\alpha
\tan\frac{\pi\alpha}{2}\beta(\mathbf{s})\left( u^{1-\alpha}-1 \right)
\right\}\bigg]\\
&& \hspace{4cm}
u^{d-1} (\sigma(\mathbf{s}))^{-d} du d \mathbf{s}.
\end{eqnarray*}
Taking the real part of the integrand, we obtain
\begin{equation*}
f_{\alpha,d}^M(\mathbf{x})=\int_{\mathbb{S}^d}
g_{\alpha,d}^M \left( \frac{ \langle \mathbf{x},\mathbf{s} \rangle
-\mu_M^M(\mathbf{s})}{\sigma(\mathbf{s})},\beta(\mathbf{s})\right)
(\sigma(\mathbf{s}))^{-d}\ d \mathbf{s}.
\end{equation*}
Note that the same conclusion holds if we let $\alpha \to 1$ in the proof due to continuity at $\alpha=1$. The direct but tedious calculation
gives the representation $\mu_M^M(\mathbf{t})$ at $\alpha=1$.
\hfill $\Box$

\bigskip\noindent
{\bf Acknowledgment:} \ The authors are grateful to Prof.\ Makoto
Maejima for very helpful comments.

\bibliographystyle{plainnat}
\bibliography{references}

\end{document}